\documentclass[a4paper,12pt]{amsart}

\usepackage[utf8]{inputenc}

\usepackage{geometry}
\usepackage{geometry}
\usepackage{amsmath,amsfonts,amsthm, amssymb,mathtools}
\usepackage{hyperref, color}

\usepackage[
backend=biber,
style=numeric-comp,
sorting=nyt,doi=false,isbn=false,url=false,eprint=false, maxbibnames=10, maxcitenames=10
]{biblatex}

\usepackage{soul}

\newcommand{\C}{\mathbb{C}}

\newcommand{\bZ}{\mathbb{Z}}

\newcommand{\bzn}{\mathbb{Z}_N}
\newcommand{\bznn}{\mathbb{Z}_N\times\mathbb{Z}_N}

\newcommand{\cL}{\mathcal{L}}
\newcommand{\cU}{\mathcal{U}}

\let\originalleft\left
\let\originalright\right
\renewcommand{\left}{\mathopen{}\mathclose\bgroup\originalleft}
\renewcommand{\right}{\aftergroup\egroup\originalright}

\newtheorem{theorem}{Theorem}[section]

\newtheorem{proposition}[theorem]{Proposition}
\newtheorem{remark}[theorem]{Remark}

\numberwithin{equation}{section}

\title[The uncertainty principle on finite cyclic groups]{The uncertainty principle for the short-time Fourier transform on finite cyclic groups: cases of equality}

\author[F. Nicola]{Fabio Nicola}
\address{Dipartimento di Scienze Matematiche, Politecnico di Torino, Corso Duca degli Abruzzi 24, 10129 Torino, Italy}
\email{fabio.nicola@polito.it}

\date{}

\begin{document}
\begin{abstract}
A well-known version of the uncertainty principle on the cyclic group $\bzn$ states that for any couple of functions $f,g\in\ell^2(\bzn)\setminus\{0\}$, the short-time Fourier transform $V_g f$ has support of cardinality at least $N$. This result can be regarded as a time-frequency version of the celebrated Donoho-Stark uncertainty principle on $\bzn$. Unlike the Donoho-Stark principle, however, 
a complete identification of the extremals is still missing. In this note we provide an answer to this problem by proving that the support of $V_g f$ has cardinality $N$ if and only if it is a coset of a subgroup of order $N$ of $\bznn$. Also, we completely identify the corresponding extremal functions $f,g$. Besides translations and modulations, the symmetries of the problem are encoded by certain metaplectic operators associated with elements of ${\rm SL}(2,\bZ_{N/a})$, where $a$ is a divisor of $N$.
Partial generalizations are given to finite Abelian groups.  
\end{abstract}

\maketitle

\section{Introduction and discussion of the main result}
For $N\geq 1$ integer, consider the cyclic group $\bzn=\bZ/N\bZ$, equipped with the counting measure.

The celebrated Donoho-Stark uncertainty principle on $\bzn$ \cite{donoho,matolscsi} (see also the seminal paper \cite{cowling}) states that, for a function $f\in \ell^2(\bzn)\setminus\{0\}$,
\begin{equation}\label{eq ds}
|{\rm supp}\, (f)|\times |{\rm supp}\, (\widehat{f})|\geq N,
\end{equation}
where ${\rm supp}\, (f)=\{j\in\bzn: f(j)\not=0\}$,  $\widehat{f}$ denotes the Fourier transform of $f$ in $\bzn$ and $|A|$ stands for the cardinality of a set $A$. Moreover, the extremal functions $f$, for which equality occurs in \eqref{eq ds}, were identified in \cite{donoho} as  characteristic functions of subgroups of $\bzn$ up to multiplication by a constant, translation and modulation (see below for the relevant definitions). 

When $N$ is prime, Terence Tao \cite{tao} then improved the above estimate as
\[
|{\rm supp}\, (f)|+ |{\rm supp}\, (\widehat{f})|\geq N+1.
\]
Refinements and generalizations to finite Abelian groups have been extensively studied  \cite{bonami,garcia,ghobber,matusiak, meshulam, smith,tao2}; see also \cite{gorjan,meshulam1992} for  extensions to compact groups and \cite{wigderson} for an up-to-date and illuminating account (and further generalizations).

Notice that the left-hand side of \eqref{eq ds} represents the cardinality of the support, in  $\bznn$, of the joint time-frequency distribution $f\otimes\widehat{f}$. Similar uncertainty inequalities hold for other time-frequency distributions, in particular for the short-time Fourier transform, which is a popular choice in signal processing \cite{mallat_book}, harmonic analysis \cite{grochenig_book} and also mathematical physics, where it is also known as coherent state transform \cite{lieb_book}. To provide its definition in the discrete setting we first introduce some notation (cf. \cite{feichtinger},\cite{krahmer}).  

For $j,k\in\bzn$ we define the translation and modulation operators $T_j$ and $M_k$ on $\ell^2(\bzn)$, and the corresponding time-frequency shifts $\pi(j,k)$ as 
\begin{equation}\label{eq transl}
T_j f(\ell)=f(\ell-j),\quad M_k f(\ell)=e^{2\pi i k\ell/N} f(\ell),\quad \pi(j,k) f= M_k T_j f,
\end{equation}
where $\ell\in\bzn$. For $f,g\in \ell^2(\bzn)$, the short-time Fourier transform of $f$ with window $g$ is the complex-valued function on $\bznn$ given by
\begin{align}\label{eq vgf}
V_g f(j,k)&= \frac{1}{\sqrt{N}}\langle f, \pi(j,k)g \rangle\\
&=\frac{1}{\sqrt{N}}\sum_{\ell\in\bzn} e^{-2\pi i k \ell/N} f(\ell) \overline{g(\ell -j)}\nonumber\qquad j,k\in\bzn.
\end{align}
The following result \cite{krahmer} is the expected counterpart of \eqref{eq ds} for the short-time Fourier transform on $\bzn$. 
\begin{theorem}\label{thm mainthm0}
If $f,g\in\ell^2(\bzn)\setminus\{0\}$ then 
\begin{equation}\label {eq 1}
|{\rm supp}\, (V_g f)|\geq N.
\end{equation}
\end{theorem}

This lower bound can be regarded as a discrete version of the so-called weak uncertainty principle \cite{groechenig} for the short-time Fourier transform in $\mathbb{R}^d$. 
It should emphasized, however, that \eqref{eq 1} does not have an exact counterpart in $\mathbb{R}^d$, because in that setting $V_g f$ cannot be  fully concentrated on a subset of finite measure of $\mathbb{R}^d\times \mathbb{R}^d$ (see e.g. \cite{grochenig_book}). In that case one considers, instead, subsets where $V_gf$ is, say, $\varepsilon$-concentrated \cite{abreu,boggiatto,bonami_demange_jaming,demange,fefferman, groechenig,lerner, lieb,romero}. Incidentally,  the corresponding sharp uncertainty principle, when the window $g$ is a Gaussian function, was  proved only recently in \cite{nicola_tilli}.

Now, despite the simplicity of the lower bound \eqref{eq 1} and its formal similarity with the Donoho-Stark inequality \eqref{eq ds}, a complete identification of the functions $f,g$ for which equality occurs in \eqref{eq 1} is still missing, except for an interesting special case studied in \cite{galbis} (see below).
In this note, partially motivated by the study of the cases of equality in the continuous-time case \cite{nicola_tilli}, we address this problem and provide a complete answer. In order to state our result, we define yet another family of unitary operators defined on certain invariant subspaces of $\ell^2(\bzn)$.  

Let $a$ be a (positive) divisor of $N$ and let $f\in\ell^2(\bzn)$ be a function supported in the subgroup $H_a:= \{ma: m=0,\ldots,N/a-1\}$ of $\bzn$ generated by $a$. For $p\in \mathbb{Z}$ we define the pointwise multiplication by a discrete ``virtual chirp" $C_{p,a}$ as 
\begin{equation}\label{eq chirp} 
C_{p,a} f(\ell)=e^{\frac{\pi i}{N}\frac{p\ell^2}{a} \big( 1+\frac{N}{a}\big)}f(\ell)\qquad \ell\in\bzn.
\end{equation}
One can easily verify that this definition makes sense (in $\bzn$), i.e. the function $C_{p,a} f$ is $N$-periodic, because of the above condition on the support of $f$ and the counterterm involving $N/a$ (whereas the exponential function alone is not well defined -- hence the name ``virtual"). The reader acquainted with the theory of metaplectic operators (see e.g. \cite{feichtinger2}) will notice that this is \textit{not} a metaplectic operator on $\bzn$ (if $a>1$), but it is unitarily equivalent to a metaplectic operator on $\mathbb{Z}_{N/a}$, via the natural identification of the linear subspace of $\ell^2(\bzn)$ of functions supported in $H_a$ with $\ell^2(\mathbb{Z}_{N/a})$; see \eqref{eq fei} below. In any case, no knowledge of the theory of metaplectic operators is needed to understand this paper.  

The following result provides the desired indentification of the subsets of $\bznn$ of smallest possibile cardinality where $V_g f$ may be fully concentrated, with corresponding extremals $f,g$.

Let $\chi_A$ denote the characteristic function of a set $A$. 

\begin{theorem}\label{thm mainthm1}
Let $f,g\in\ell^2(\bzn)$. The following statements are equivalent. 
\begin{itemize}
    \item[(a)] $|{\rm supp}\, (V_g f)|=N$.
    \item[(b)] ${\rm supp}\, (V_g f)$ is a coset of a subgroup of $\bznn$ of order $N$.
    \item[(c)] There exist a divisor $b$ of $N$,  $p\in\{0,\ldots,b-1\}$, $\lambda,\mu\in\bznn$, $c_1,c_2\in\mathbb{C}\setminus\{0\}$ such that 
    \[
    g= c_1 \pi(\lambda) C_{p,a} \chi_{H_a},\quad f=c_2 \pi(\mu) g,
    \]
    where $a=N/b$ and $H_a\subseteq \bzn$ is the subgroup generated by $a$. 
\end{itemize}
\end{theorem}
It will follow from the proof that, with $f,g$ as in Theorem \ref{thm mainthm1} (c), 
\[
{\rm supp}\, (V_g f)=\mu+H_{b,p},
\]
where (with $a=N/b$, as above)
\begin{equation}\label{eq def hbp}
H_{b,p}:=\{(ma,nb+mp):\ m=0,\ldots,b-1,\ n=0,\ldots,a-1\},
\end{equation}
i.e., the lattice generated by $(a,p)$ and $(0,b)$ in $\bznn$. 
Indeed, it is known that all the subgroups of order $N$ of $\bznn$ have this form for some divisor $b$ of $N$ and $p\in\{0,\ldots,b-1\}$ (see e.g. \cite{hampejs}).

We will see that the operator $C_{p,a}$, $p\in\{0,\ldots,N/a-1\}$, is associated with the phase-space map $H_{b,0}\to H_{b,p}$, $(ma,nb)\mapsto(ma,nb+mp)$ ($m=0,\ldots,b-1,\ n=0,\ldots,a-1$), in the sense that the intertwining properties \eqref{eq intet} and \eqref{eq intet2} below hold true. Notice that, however, this map is not the restriction to $H_{b,0}$ of an element of ${\rm SL}(2,\bzn)$ ($a$ is not supposed to be a divisor of $p$ in $\bzn$), nor it is a group  isomorphism (it is just a bijection). Indeed, as a metaplectic operator in $\mathbb{Z}_{N/a}$ (via the above mentioned identification), $C_{p,a}$ is associated with the matrix $\begin{pmatrix}
1&0\\
p&1
\end{pmatrix}\in{\rm SL}(2,\mathbb{Z}_{N/a})$ (cf. \cite{feichtinger2} and \eqref{eq fei} below).

 The extremal functions $f,g$ satisfying the additional condition $|{\rm supp}\, (f)|+ |{\rm supp}\, (g)|\geq N+1$ 
 were already identified in \cite[Theorem 3]{galbis}, and turn out to be those in Theorem \ref{thm mainthm1} (c) for $b=N$, therefore $a=1$ (hence, in that case the multiplication by a true chirp function appears -- as opposed to the genuinely virtual chirps that we alluded to above).
 The proof of Theorem \ref{thm mainthm1} in full generality however will require a different approach, to take into account more tightly of the symmetries of the problem -- which correspond to the operators appearing in Theorem \ref{thm mainthm1} (c).   
 
It is known that Theorem \ref{thm mainthm0} generalizes to finite Abelian groups (see \cite[Proposition 4.1]{krahmer} and Section \ref{sec conc} below). The cases of equality seem certainly worthy studying in that framework too, as well as for the other discrete uncertainty principles appearing in \cite{krahmer}.  While we postpone this investigation to a future work, in Section \ref{sec conc} we briefly show that the implication (a)$\Longrightarrow$(b) in Theorem \ref{thm mainthm1} indeed generalizes (easily) to arbitrary finite Abelian groups. 

We conclude by observing that subgroups (or structured subsets) $S\subset \bznn$ play also an important role in the construction of orthonormal basis of $\ell^2(\bzn)$ of the type $\{\pi(j,k) g:\ (j,k)\in S$\}, for some $g\in\ell^2(\bzn)$. This problem has a long tradition, for which we address to the recent contribution \cite{zhou} and the references therein. Also, the study of uncertainty principles involving the cardinality of the support of discrete  time-frequency distributions has interesting applications  in signal recovery theory, for which we address to \cite{krahmer,tao2}. 
 
Briefly, this note is organized as follows. In Section \ref{sec prel} we recall some basic formulas from time-frequency analysis on $\bzn$, which will be needed to fully exploit the symmetries of the problem. Section \ref{sec proof} is devoted to the proof of Theorem \ref{thm mainthm1}. Finally in Section \ref{sec conc} we report on the very short proof of Theorem \ref{thm mainthm0}, for the sake of completeness, and we discuss the above mentioned partial generalisation of Theorem \ref{thm mainthm1} to finite Abelian groups.  

\section{Notation and preliminary results}\label{sec prel}
\subsection{Notation}
We denote by $\bzn=\mathbb{Z}/N\mathbb{Z}$ the cyclic group of order $N$, equipped with the counting measure. The inner product and corresponding norm in $\ell^2(\bzn)\equiv \mathbb{C}^N$ are denoted by $\langle\cdot,\cdot\rangle$ and $\|\cdot\|_{\ell^2(\bzn)}$ respectively. The support of a function $f\in\ell^2(\bzn)$ is denoted by ${\rm supp}\, (f)=\{j\in\bzn: f(j)\not=0\}$, and similarly for functions on $\bznn$. We denote by $|A|$ the cardinality of a set $A$, and by $\chi_A$ its characteristic function. 

In the following $b$ will always denote a (positive) divisor of $N$ and $a=N/b$. $H_a$ stands for the subgroup of $\bzn$ generated by $a$, whereas $H_{b,p}$, $p\in\{0,\ldots,b-1\}$, is the subgroup of order $N$ of $\bznn$ defined in \eqref{eq def hbp}. 

We already defined in \eqref{eq transl} the translation operators $T_j$, $j\in\bzn$, the modulation operators $M_k$, $k\in\bzn$ and the time-frequency shifts $\pi(j,k)= M_k T_j$, as unitary operators on $\ell^2(\bzn)$. We also defined the chirp operator $C_{p,a}$ in \eqref{eq chirp}, as a unitary operator defined on the linear subspace of $\ell^2(\bzn)$ of functions supported in $H_a$. The short-time Fourier transform $V_g f$, for $f,g\in\ell^2(\bzn)$, was defined in \eqref{eq vgf}.  

\subsection{Preliminaries from time-frequency analysis on $\bzn$}
We collect here some elementary formulas that will be useful in the following. We only sketch the proofs, or we omit them completely whenever they are straightforward computations (for analogous results in $\mathbb{R}^d$, see \cite{grochenig_book}). 

In the following, $f,g$ denote functions in $\ell^2(\bzn)$. 

First of all, we observe that the time-frequency shifts enjoy the commutation relations
\begin{equation}\label{eq proj}
\pi(j,k)\pi(j',k')=e^{2\pi i(kj'-k' j)/N}\pi(j',k')\pi(j,k)
\end{equation}
for $j,k,j',k'\in\bzn$.

The short-time Fourier transform $V_g f(j,\cdot)$ can be regarded as the Fourier transform on $\bzn$ of $f \overline{T_j g}$. Hence, applying the Plancherel theorem for the Fourier transform, followed by the Fubini Theorem, one obtains at once the Parseval equality for the short-time Fourier transform, which reads
\begin{equation}\label{eq energy}
\|V_g f\|_{\ell^2(\bznn)}=\|f\|_{\ell^2(\bzn)} \|g\|_{\ell^2(\bzn)}.
\end{equation}
We also have the following pointwise estimate
\begin{equation}\label{eq cs}
|V_g f(j,k)|\leq \frac{1}{\sqrt{N}}\|f\|_{\ell^2(\bzn)} \|g\|_{\ell^2(\bzn)}\qquad j,k\in\bzn,
\end{equation}
which is an immediate consequence of the Cauchy-Schwarz inequality.

The following covariance-type properties can be checked by direct computation, using the definition of $V_gf$ and \eqref{eq proj}:  
\begin{equation}\label{eq covtf}
V_{g}(\pi(j,k)f)(j',k')= e^{2\pi i (k-k')j/N} V_g f(j'-j,k'-k)
\end{equation}
and
\begin{equation}\label{eq covamb}
V_{\pi(j,k)g}(\pi(j,k)f)(j',k')= e^{2\pi i(kj'-k'j)/N} V_g(f)(j',k'),
\end{equation}
for $j,k,j',k'\in\bzn$. 

Similarly, if $ N=ab$, we have the  intertwining property (cf. \eqref{eq chirp})
\begin{equation}\label{eq intet}
C_{-p,a}\, \pi(ma,nb+mp)= e^{\pi i pm^2(1+b)/b}
\pi(ma,nb)\, C_{-p,a}
\end{equation}
for $p\in\{0,\ldots,b-1\}$, $m\in\{0,\ldots,b-1\}$, $n\in\{0,\ldots,a-1\}$. Equivalently,
\begin{equation}\label{eq intet2}
\pi(ma,nb+mp)\,C_{p,a}=e^{\pi i pm^2(1+b)/b}  C_{p,a}\, \pi(ma,nb).
\end{equation}
Unlike the previous formulas, \eqref{eq intet} and \eqref{eq intet2} hold, however, only on the subspace of $\ell^2(\bzn)$ of functions supported in $H_a$, which indeed is an invariant subspace for all the operators appearing in \eqref{eq intet} and \eqref{eq intet2} (when $a=1$, hence $b=N$, $C_{p,a}f$ is defined for every $f\in\ell^2(\bzn)$ -- in that case $C_{p,a}$ is a metaplectic operator on $\bzn$ \cite{feichtinger2}, and \eqref{eq intet2} reduces to \cite[Lemma 3.1 (iii)]{feichtinger2}). 

\begin{remark}
More generally, let $\cL_a\subseteq\ell^2(\bzn)$ be the subspace of functions supported in $H_a$, and let $\cU_a:\cL_a\to\ell^2(\bZ_b)$ be given by $\cU_a f(m)=f(ma)$, $m\in\bZ_b$ ($b=N/a$). Then 
\begin{equation}\label{eq fei}
    \cU_a C_{p,a}\cU_a^{-1}f(m)=e^{
    \pi i p m^2(1+b)/b} f(m)\qquad m\in\bZ_b,
\end{equation}
is a metaplectic operator on $\bZ_b$, as defined in \cite[Section 3 (iii)]{feichtinger2}. 
\end{remark}

\begin{remark}\label{rem inversion}
It is easy to see that $f\in\ell^2(\bzn)$ is determined by $V_f f$ up to multiplication by a complex number (of modulus 1 if $\|f\|_{\ell^2(\bzn)}$ is given, by \eqref{eq energy}); cf. \cite[Section 4.2]{grochenig_book}. Indeed, from the definition of $V_g f$, using the inversion formula for the discrete Fourier transform we obtain
\[
\frac{1}{\sqrt{N}}\sum_{k\in\bzn} e^{2\pi i k j'/N} V_f f(j,k)= f(j') \overline{f(j'-j)}\qquad j,j'\in\bzn.
\]
This already tells us that the support of $f$ is determined from the function $V_f f$. Moreover, if   $f(j_0)\not=0$, say,  choosing $j'=j_0+j$ yields
\[
f(j_0+j)= \frac{1}{\overline{f(j_0)}\sqrt{N}}\sum_{k\in\bzn} e^{2\pi i k (j_0+j)/N} V_f f(j,k)\qquad j\in\bzn,
\]
which gives the claim. 
\end{remark}

\section{Proof of the main result (Theorem \ref{thm mainthm1})}\label{sec proof}
We can suppose that $f$ and $g$ are normalized in $\ell^2$; hence $\|f\|_{\ell^2(\bzn)}=\|g\|_{\ell^2(\bzn)}=1$.
\bigskip

(a) $\Longrightarrow$ (b) 

It follows from \eqref{eq energy} and \eqref{eq cs} that, if $|{\rm supp}\, (V_g f)|=N$, then 
\[
|\langle f, \pi(\lambda)g\rangle|=\sqrt{N}|V_g f(\lambda)|=1\qquad \lambda\in {\rm supp}\, (V_g f).
\]
Hence,
\begin{equation}\label{eq 3}
f=c(\lambda) \pi(\lambda) g\qquad \lambda\in {\rm supp}\, (V_g f)
\end{equation}
for some $c(\lambda)\in\mathbb{C}$, $|c(\lambda)|=1$. 
Applying the short-time Fourier transform $V_g$ to both sides of \eqref{eq 3}  and using \eqref{eq covtf} we obtain
\[
|V_g f(\mu)|= | V_g g(\mu-\lambda)|\qquad \lambda\in {\rm supp}\, (V_g f),\ \mu\in\bzn.
\]
Hence 
\begin{equation}\label{eq 2}
{\rm supp}\, (V_g f)=\lambda+{\rm supp}\, (V_g g)\qquad \lambda\in {\rm supp}\, (V_g f). 
\end{equation}
We claim that $H:={\rm supp}\, (V_g g)$ is a subgroup of $\bznn$. Indeed, let $\mu_1,\mu_2\in H$. For any $\lambda\in {\rm supp}\, (V_g f)$ we can write $\mu_1=\lambda_1-\lambda$ and $\mu_2=\lambda_2-\lambda$ for some $\lambda_1,\lambda_2\in {\rm supp}\, (V_g f)$. Then $\mu_1-\mu_2=\lambda_1-\lambda_2$ belongs to $H$ thanks to \eqref{eq 2} (applied with $\lambda_2$ in place of $\lambda$). 

\bigskip
(b) $\Longrightarrow$ (c) 
\nopagebreak

Suppose that ${\rm supp}\, (V_g f)=\lambda_0+H$, where $H$ is a subgroup of $\bznn$ of order $N$ and $\lambda_0\in\bznn$. Then $|{\rm supp}\, (V_g f)|=N$ and in particular we can take for granted what we have just proved. Hence, since $\lambda_0\in {\rm supp}\, (V_g f)$, from \eqref{eq 2} we deduce that $H={\rm supp}(V_g g)$. From the classification of the subgroups of $\bznn$ (see e.g. \cite[Theorem 1]{hampejs}) we see that there exist a divisor $b$ of $N$ and $p\in\{0,\ldots,b-1\}$ such that $H=H_{b,p}$ as in \eqref{eq def hbp}, where $a=N/b$. 

By \eqref{eq 3} we  have
\[
\pi(\lambda) g =c(\lambda,\lambda') \pi(\lambda') g\qquad \lambda,\lambda'\in \lambda_0+H_{b,p}
\]
for some factor $|c(\lambda,\lambda')|=1$, or equivalently, using \eqref{eq proj},
\begin{equation}\label{eq 5}
\pi(\lambda)g=c(\lambda)g \qquad \lambda\in H_{b,p}
\end{equation}
with $|c(\lambda)|=1$. 

Applying this formula repeatedly with $\lambda=(0,b)$ we see that 
\[
M_{nb} g = c^n g
\qquad n\in\mathbb{Z}
\]
for some constant $c\in\mathbb{C}$ with $c^a=1$, hence $c=e^{2\pi i j_0/a}=e^{2\pi i bj_0/N}$ for some $j_0\in \{0,\ldots,a-1\}$. In particular for $n=1$ we obtain
\[
e^{2\pi i b\ell/N} g(\ell)= e^{2\pi i bj_0/N} g(\ell) \qquad \ell\in\bzn,
\]
and therefore $g$ is supported in the coset $\{j_0+ma:\ m=0,\ldots ,b-1\}$. 

To reduce things to the case of a product-type subgroup of $\bznn$, we introduce the function
\begin{equation}\label{eq gammag}
\gamma:= C_{-p,a} T_{-j_0} g.
\end{equation}
Notice that $T_{-j_0} g$, and therefore $\gamma$, is supported in the subgroup $H_a$ generated by $a$ in $\bzn$. 

Applying the operator $C_{-p,a}T_{-j_0}$ to both sides of \eqref{eq 5} and using \eqref{eq proj} and \eqref{eq intet} we obtain  
\[
\pi(\lambda)\gamma= c(\lambda) \gamma \qquad \lambda \in H_{b,0}
\]
for a new constant $c(\lambda)\in\mathbb{C}$, $|c(\lambda)|=1$.

Applying the short-Fourier transform $V_\gamma$ to both sides of this equality we obtain
\[
V_\gamma \gamma (\mu-\lambda)=  V_\gamma (\pi(\lambda)\gamma)(\mu)= c(\lambda) V_\gamma\gamma(\mu) \qquad \lambda,\mu\in H_{b,0},
\]
for the same constant $c(\lambda)$ (which is therefore independent of $\mu)$, where the first equality follows from \eqref{eq covtf} (the exponential factor appearing there is $=1$ because $\lambda,\mu\in H_{b,0}$).

Applying repeatedly this formula with $\lambda=(a,0)$ or $\lambda=(0,b)$ we obtain, for $m,n\in\mathbb{Z}$, 
\begin{equation}\label{eq 6}
V_\gamma \gamma (-ma,-nb)= c(a,0)^m c(0,b)^nV_\gamma\gamma(0)=c(a,0)^m c(0,b)^n/\sqrt{N},
\end{equation}
because $V_\gamma\gamma(0)=\|\gamma\|^2_{\ell^2(\bzn)}/\sqrt{N}=1/\sqrt{N}$. Moreover $c(a,0)^b=1$ and $c(0,b)^a=1$, namely there exist $j_1\in\{0,\ldots,a-1\}$ and $k_1\in \{0,\ldots,b-1\}$ such that
\begin{equation}\label{eq const}
c(a,0)=e^{2\pi i k_1/b}=e^{2\pi i k_1a/N},\quad c(0,b)=e^{2\pi i j_1/a}=e^{2\pi i j_1b/N}. 
\end{equation}
This also tells us that $V_\gamma \gamma(\lambda)=0$ for $\lambda\not\in H_{b,0}$, because $|V_\gamma \gamma|^2=1/N$ on $H_{b,0}$ (by \eqref{eq 6} and \eqref{eq const}), $|H_{b,0}|=ab=N$ and  $\|V_\gamma \gamma\|_{\ell^2(\bznn)} =\|\gamma\|^2_{\ell^2(\bznn)}=1$ by \eqref{eq energy}. Hence $
{\rm supp}\, (V_\gamma\gamma)= H_{b,0}$.

Summing up, 
\begin{equation}\label{eq 12}
    V_{\gamma}\gamma(j,k)=\frac{1}{\sqrt{N}}e^{-2\pi i(k_1 j+k j_1)/N} \chi_{H_{b,0}}(j,k)
\end{equation}
for $(j,k)\in\bznn$.

Now, to identify the function $\gamma$ (up to a phase factor) we are going to exhibit a function $\tilde{\gamma}$ such that $V_{\tilde{\gamma}}\tilde{\gamma}=V_\gamma\gamma$; this will give $\tilde{\gamma}=c\gamma$ for some $|c|=1$ by Remark \ref{rem inversion}. 

To this end, consider the subgroup $H_a$ of $\bzn$ generated by $a$. An explicit computation shows that 
\[
\frac{1}{b}V_{\chi_{H_a}}\chi_{H_a}=\frac{1}{\sqrt{N}}\chi_{H_{b,0}}.
\]
Hence setting
\[
\tilde{\gamma}= \frac{1}{\sqrt{b}} M_{-k_1}T_{j_1} \chi_{H_a} 
\]
and using \eqref{eq covamb} we obtain 
\[
V_{\tilde{\gamma}}\tilde{\gamma}(j,k)=\frac{1}{\sqrt{N}}e^{-2\pi i(k_1 j+k j_1)/N} \chi_{H_{b,0}}(j,k)
\]
for $(j,k)\in\bznn$. Hence $V_{\tilde{\gamma}}\tilde{\gamma}=V_\gamma\gamma$ by \eqref{eq 12}. 

By Remark \ref{rem inversion} we deduce that 
\[
\gamma=\frac{c}{\sqrt{b}} M_{-k_1}T_{j_1} \chi_{H_a}
\]
for some constant $c\in\C$, $|c|=1$.
Moreover, since $\gamma$ is supported in $H_a$ we have in fact $j_1=0$. Coming back to the function $g$ (cf. \eqref{eq gammag}) we obtain 
\[
C_{-p,a}T_{-j_0} g=\frac{c}{\sqrt{b}} M_{-k_1} \chi_{H_a}
\]
and therefore, since $C_{p,a}$ and $M_{k_1}$ commute and using \eqref{eq proj}, 
\[
g=\frac{c}{\sqrt{b}} M_{-k_1} T_{j_0} C_{p,a}\chi_{H_a}
\]
for a new constant $c\in\C$, $|c|=1$.

By \eqref{eq 3}, also $f$ has the desired form. 

\bigskip
(c) $\Longrightarrow$ (a) 
\nopagebreak

Let $f,g$ be as in the statement (c). The computation of the support of $V_g f$ is straightforward. Indeed, let $
h\coloneqq\chi_{H_a}/\sqrt{b}$.
We have already observed that  $V_{h}h=\chi_{H_{b,0}}/\sqrt{N}$. Moreover, for  $m\in\{0,\ldots,b-1\}$, $n\in\{0,\ldots, a-1\}$ we have 
\begin{align*}
V_{C_{p,a}h} C_{p,a}h(ma,nb+mp)&=\langle C_{p,a}h,\pi(ma,nb+mp)C_{p,a}h\rangle/\sqrt{N} \\
&=c\,\langle C_{p,a}h,C_{p,a}\pi(ma,nb)h\rangle/\sqrt{N} \\
&=c\, V_{h}h (ma,nb) =c/\sqrt{N}
\end{align*}
for some constant $c$, $|c|=1$, where we used \eqref{eq intet2} and the fact that $C_{p,a}$ is unitary on the subspace of functions supported in $H_a$.

Since, by \eqref{eq energy},
\[
\|V_{C_{p,a}h} C_{p,a}h\|_{\ell^2(\bznn)}=\|C_{p,a}h\|^2_{\ell^2(\bzn)}=\|h\|^2_{\ell^2(\bzn)}=1,
\]
we have ${\rm supp}\, (V_{C_{p,a}h} C_{p,a}h)=H_{b,p}$. 

Using \eqref{eq covamb} we deduce that 
\[
{\rm supp}\, (V_{g} g)=H_{b,p}
\]
and by \eqref{eq covtf}
\[
{\rm supp}\, (V_{g}f)=\mu+{\rm supp}\, (V_{g}g) =\mu+H_{b,p},
\]
which has therefore cardinality $N$. 

This concludes the proof of Theorem \ref{thm mainthm1}.

\section{Concluding remarks}\label{sec conc}
In this section we report on the short proof of Theorem \ref{thm mainthm0}, following \cite[Proposition 4.1]{krahmer}. We also discuss a partial generalization of Theorem \ref{thm mainthm1} to the case of arbitrary finite Abelian groups.  


\subsection{Proof of Theorem \ref{thm mainthm0}}(\cite[Proposition 4.1]{krahmer}) 
It follows from  \eqref{eq cs} and \eqref{eq energy} that, if $S={\rm supp\,} (V_g f)$, 
\[
\|f\|^2_{\ell^2(\bzn)}\|g\|^2_{\ell^2(\bzn)}
= \sum_{(j,k)\in S}|V_g f(j,k)|^2\leq 
\frac{|S|}{N}\|f\|^2_{\ell^2(\bzn)}\|g\|^2_{\ell^2(\bzn)}
\]
so that $|S|\geq N$ if $f,g\not=0$.  

\subsection{Generalization to finite Abelian groups} The basic definitions from time-frequency analysis on $\mathbb{R}^d$ have a natural counterpart on finite Abelian groups (see e.g. \cite{feichtinger,krahmer}). 

In short, on a finite Abelian group $A$ (equipped with the counting measure) one can define the translation operator $T_j$, $j\in A$, and modulation operator $M_k$, $k\in \widehat{A}$ (the dual group), as well as the corresponding time-frequency shifts $\pi(j,k)$ as unitary operators on $\ell^2(A)$ by 
\[
T_j f(\ell)=f(\ell-j),\quad M_k f(\ell)=e^{2\pi i \langle k,\ell\rangle} f(\ell),\quad \pi(j,k) f= M_k T_j f,\quad \ell\in A,
\]
where the map $A\ni\ell \mapsto\langle k,\ell\rangle$ denotes the {\it additive} character $k\in\widehat{A}$ (that is a continuous  homomorphism $A\to \mathbb{R}/\mathbb{Z}$) -- so that  $e^{2\pi i \langle k,\ell\rangle}$ is the corresponding multiplicative character (a continuous  homomorphism $A\to \mathbb{T}=\{z\in\mathbb{C}:\ |z|=1\}$ -- the circle group). 

We then define the short-time Fourier transform as
\[
V_g f(j,k)=\frac{1}{\sqrt{|A|}}\langle f,M_k T_j g\rangle_{\ell^2(A)},
\]
for $g,f\in\ell^2(A)$, $(j,k)\in A\times\widehat{A}$, and the properties \eqref{eq energy} and \eqref{eq cs} generalize in the obvious way (cf. \cite{feichtinger,krahmer}). Hence, also the above proof of Theorem \ref{thm mainthm0} extends to the case of finite Abelian groups, as already observed in \cite[Proposition 4.1]{krahmer}.

Similarly, one can easily check that the formulas \eqref{eq proj} and \eqref{eq covtf} generalize naturally, and also the proof of (a)$\Longrightarrow$(b) in Theorem \ref{thm mainthm1} carries on essentially without changes. Summing up, we have the following result. 
\begin{proposition}\label{pro 1}
Let $A$ be a finite Abelian group, and $f,g\in\ell^2(A)\setminus\{0\}$. Then $|{\rm supp}\, (V_g f)|\geq |A|$. If $|{\rm supp}\, (V_g f)|=|A|$ then ${\rm supp}\, (V_g f)$ is a coset of a subgroup of $A\times \widehat{A}$ of order $N$. 
\end{proposition}

\section*{Acknowledgments} 
It is a pleasure to thank Paolo Tilli and S. Ivan Trapasso for interesting discussions about the subject of this paper.


\printbibliography
\end{document}